\documentclass{article}

\usepackage{amsmath}
\usepackage{amssymb}
\usepackage{graphicx}
\usepackage{epic,eepic}

\usepackage{color}

\usepackage{amstext,amsthm,amssymb,amsmath}
\usepackage{epsfig,subfigure,epstopdf}
\usepackage{graphicx}
\usepackage{subfigure}
\usepackage{wrapfig}
\usepackage{fullpage}
\usepackage{color}
\usepackage{multirow}
\usepackage{tabulary}
\usepackage{booktabs}
\usepackage{enumerate}

\newtheorem{theorem}{Theorem}
\newtheorem{lemma}{Lemma}
\newtheorem{remark}{Remark}

\usepackage[utf8]{inputenc}
\usepackage[english]{babel}
\usepackage{hyperref}

\usepackage{nomencl}

\theoremstyle{definition}

\begin{document}
\title{Fast Online Generalized Multiscale Finite Element Method using Constraint Energy Minimization}

\author{
Eric T. Chung \thanks{Department of Mathematics,
The Chinese University of Hong Kong (CUHK), Hong Kong SAR. Email: {\tt tschung@math.cuhk.edu.hk}.
The research of Eric Chung is supported by Hong Kong RGC General Research Fund (Project 14317516).}
\and
Yalchin Efendiev \thanks{Department of Mathematics \& Institute for Scientific Computation (ISC),
Texas A\&M University,
College Station, Texas, USA. Email: {\tt efendiev@math.tamu.edu}.}
\and
Wing Tat Leung \thanks{Department of Mathematics, Texas A\&M University, College Station, TX 77843}
}

\maketitle

\begin{abstract}
Local multiscale methods often construct multiscale basis functions
in the offline stage without taking into account input parameters,
such as source terms, boundary conditions, and so on. These basis
functions are then used in the online stage with a specific
input parameter to solve the global problem at a reduced
computational cost.
Recently, online approaches have been introduced, where
multiscale basis functions are adaptively constructed in
some regions to reduce the error significantly.
In multiscale methods, it is desired to have only 1-2 iterations
to reduce the error to a desired threshold.
Using Generalized Multiscale Finite Element Framework \cite{egh12}, it was shown
that by choosing
sufficient number of offline basis functions, the error
reduction can be made independent of physical parameters, such as scales
and contrast.
In this paper, our goal is to improve this. Using our recently
proposed approach \cite{chung2017constraint} and special online basis
construction in oversampled regions, we show that the error reduction
can be made sufficiently large by appropriately
selecting oversampling regions. Our numerical results show that one can achieve
a three order of magnitude error reduction, which is better
than our previous methods.
We also develop an adaptive algorithm and enrich in selected
regions with large residuals.
In our adaptive method, we show that the convergence rate
can be determined by a user-defined parameter
and we confirm this by numerical simulations.
The analysis of the method is presented.

\end{abstract}

\section{Introduction}

Many multiscale problems are prohibitively expensive to solve due to
scale disparity  and high contrast. These problems are often
solved using some type of reduced-order
models. These
include numerical homogenization \cite{dur91,weh02},
multiscale finite element methods \cite{hw97,eh09,ehw99},
heterogeneous multiscale methods \cite{ee03}, variational
multiscale methods \cite{hughes1995multiscale},
mortar multiscale methods \cite{Wheeler_mortar_MS_12},
localized orthogonal decomposition \cite{maalqvist2014localization},
and so on. The main idea behind local reduced-order model reduction
techniques is to compute multiscale basis functions in each coarse block.
These basis functions are computed using solutions of local
problems.

Various approaches have been developed for designing multiscale
basis functions. One of the earlier works
\cite{eh09,hw97,ehw99, chung2016adaptive}
 use
harmonic extensions of standard finite element basis
functions in computing multiscale basis
functions. Because of ``homogeneous'' traces along coarse boundaries,
these approaches can have large errors due to the mismatch between
the fine-grid solution and multiscale solutions along the edges of
coarse blocks. These approaches have been generalized
by using oversampling ideas \cite{hw97,ehw99},
where one uses larger regions
and solve local problems. The solutions of these local problems
are then used in constructing boundary conditions
for multiscale basis functions.
These approaches reduce the errors due to boundary conditions.

In later works \cite{ge09_2, egh12,chung2016adaptive},
the authors showed that in the presence
of high-contrast, one needs multiple basis functions.
In \cite{egh12,chung2016adaptive},
Generalized Multiscale Finite Element is introduced,
where the authors propose a systematic way of computing multiscale
basis functions. The multiscale basis functions are computed
by solving spectral problems in each coarse patch and selecting the
dominant eigenvectors. In particular, the eigenvalues are
ordered in increasing order and the eigenvectors corresponding
to small eigenvalues are selected. The spectral convergence rate
$1/\Lambda$ has been derived for these approaches, where
$\Lambda$ is the smallest eigenvalue (across all coarse blocks)
 whose corresponding eigenvector
is not included in the coarse space.
In \cite{chung2017constraint}, using oversampling ideas and localization
ideas \cite{owhadi2014polyharmonic, maalqvist2014localization, owhadi2017multigrid,hou2017},
the authors propose an approach which provides both mesh-size
dependent convergence and spectral convergence. The main idea
of this approach, called CEM-GMsFEM, is to (1) compute some GMsFEM
basis (2) use constrained energy minimization in oversampling domains
to construct multiscale basis functions. As a result,
we have a minimal number of basis functions and can show $H/\Lambda^{1/2}$
convergence rate.

The above approaches can be classified as {\it offline} methods because
the construction of multiscale basis functions does not take
into account the right hand side. The offline methods can be tuned
in various ways to achieve smaller errors; however, the error
decay slows down as we add basis functions after a certain number
of basis functions are selected. This slow down is due to some
slow decay after certain eigenvalue. To improve this, in \cite{chung2015residual,chung2015online,chan2016adaptive,chung2016onlineperforated},
the authors propose an online approach. The main idea of online approaches
is to add multiscale basis functions using the residual information
after computing the coarse-grid solution. These online multiscale basis
functions are computed adaptively and are chosen to decrease the error
the most. They are solutions of local problems.
Our analysis in \cite{chung2015residual,chung2015online}
shows that the error decay is proportional to $1-C \Lambda$,
where $C$ is the constant (independent of scales and contrast)
 that guarantees the positivity of this quantity.
 This indicates
that the error is not reduced unless $\Lambda$ is sufficiently away
from $0$, i.e., we have suffcient number of offline basis functions.
This was demonstrated analytically and numerically
in our papers \cite{chung2015residual,chung2015online}.
Since the online procedure can be
costly, our goal is to perform only 1-2 iterations.

In this paper, we would like to investigate online approaches for CEM-GMsFEM
and show that one can significantly improve the existing online approaches
for some cases. In the paper, first we present an online approach, which
differs from our previous approach
since CEM-GMsFEM uses oversampling. In particular,
the online basis functions are formulated in the oversampled regions.
Secondly, we present an analysis of the proposed method.
Our analysis shows that the error decay by adding online basis
functions can be significantly
better compared to $1-\Lambda$ in online GMsFEM. The error decay
can be made close to $0$
(i.e., we obtain very accurate approximation in one iteration)
by choosing larger oversampling regions
{\it provided} we have sufficient number of offline basis functions.
To our best knowledge, this is a first result of this kind.
Moreover, the online approaches can be made adaptive
and adaptive error indicators can be derived.

We present numerical result. In our numerical results, we consider
high-contrast permeability fields  and place the source
term in different locations. All results show that the error drops
3 orders of magnitude, which is much better compared to previous
online GMsFEM. We also present numerical results using adaptivity,
which shows that by selecting only some (few) regions, one can achieve
a significant error decay.
Moreover, our adaptive algorithm allows one to input a parameter
which specifies a desired convergence rate.

The paper is organized as follows. In Section \ref{sec:prelim},
we present some preliminaries. In Section \ref{sec:method},
we present the construction of offline basis functions,
and in Section \ref{sec:online},
we present our online adaptive enrichment algorithm.
Section \ref{sec:num} is devoted to
numerical results.
The analysis of our method is presented in Section \ref{sec:analysis}.
Finally, we present some concluding remarks
in Section \ref{sec:conclusion}.

\section{Preliminaries}
\label{sec:prelim}

In this paper, we consider a class of multiscale problems of the form
\begin{equation}
-\mbox{div}\big(\kappa(x)\,\nabla u\big)=f, \quad\text{in}\;\Omega,\label{eq:original}
\end{equation}
subject to the homogeneous Dirichlet boundary condition $u=0$ on
$\partial\Omega$, where $\Omega \subset \mathbb{R}^d$ is the computational domain. We
assume that $\kappa(x)$ is a heterogeneous coefficient with multiple
scales and very high contrast.
In solving (\ref{eq:original}), it is desirable to construct multiscale basis functions
that can be computed locally and give coarse-mesh convergence rate independent of the heterogenities and contrast.
In \cite{chung2017constraint}, such approach has been developed.
When the source term $f$ does not belong to $L^2(\Omega)$ and when one needs to obtain solutions
with more refined accuracy, it is desirable to construct online basis functions
that capture properties unresolvable by offline basis functions.
It is the purpose of this paper to do this.

\begin{figure}[ht!]
\centering
\includegraphics[width=3in, height=3in]{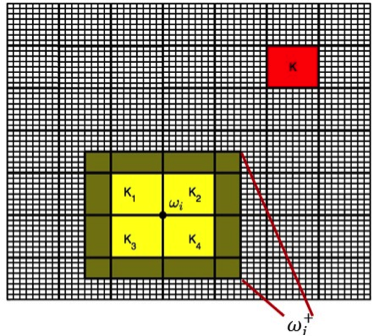}
\caption{Illustration of the coarse grid, fine grid and oversampling domain.}
\label{fig:illustration}
\end{figure}

Next, we introduce the notion of fine and
coarse grids. We let $\mathcal{T}^{H}$ be a usual conforming partition
of the computational domain $\Omega$ into $N$ finite elements (triangles,
quadrilaterals, tetrahedra, etc.), and let $H$ be the mesh size of $\mathcal{T}^H$. We refer to this partition as the
coarse grid and assume that each coarse element is partitioned into
a connected union of fine grid blocks. The fine grid partition will
be denoted by $\mathcal{T}^{h}$, and is by definition a refinement
of the coarse grid $\mathcal{T}^{H}$. Here, we use $h$ to denote the fine mesh size of $\mathcal{T}^h $.
In Figure \ref{fig:illustration}, we give an illustration of the fine grid, coarse grid, and
oversampling domain.
In the figure, the coarse grid is contained by a union of rectangular coarse elements, denoted generically by $K$.
Each coarse element is a union of finer rectangular elements.
Moreover, for each coarse grid node $x_i$, we define $\omega_i$ as the union of coarse elements having the vertex $x_i$.
We also define $\omega_i^+$ as an oversampled region for $\omega_i$.
Finally, we define $N_c$ as the number of coarse grid vertices.

We let $V=H_0^1(\Omega)$. The weak solution $u\in V$ of the problem \eqref{eq:original} satisfies
\begin{equation}
a(u,v) = (f,v), \;\;\forall v \in V,
\label{eq:global_pro}
\end{equation}
where $a(u,v)=\int_{\Omega}\kappa \nabla u\cdot \nabla v$ and $(f,v) = \int_{\Omega}fv$.
We would like to find a multiscale solution $u_{ms}$ in a subspace of $V$, denoted as $V_{ms}$, satisfying
\begin{equation}
a(u_{ms},v) = (f,v), \;\;\forall v \in V_{ms}.
\label{eq:global_ms}
\end{equation}
We will introduce the construction of the multiscale finite element space $V_{ms}$ in the next section.
We remark that the multiscale space $V_{ms}$ consists of two components, which are the offline part and the online part.
For the offline part, we will construct multiscale basis functions in the offline stage,
that is, before solving the problem (\ref{eq:global_ms}). Note that these basis functions are independent of the source term $f$.
For the online part, we will construct multiscale basis functions in the online stage
using the residual of an approximate solution. Note that, these basis functions depend on the source term $f$.
In Section \ref{sec:method}, we present the construction of offline basis functions.
In Section \ref{sec:online}, we present the construction of online basis functions
and an online adaptive enrichment algorithm.

\section{Offline basis functions}
\label{sec:method}
In this section, we present the construction of the offline multiscale finite element space. The construction of our offline basis functions follows the framework in \cite{chung2017constraint}. To construct the multiscale space, we will first construct the auxiliary space by solving local spectral problem for each coarse element $K$. Next, we will construct the multiscale basis functions by solving some local minimization problems using the auxiliary basis functions. Our multiscale finite element space will then be the span of these multiscale basis functions. Next, we will discuss the construction of both auxiliary space and multiscale space in detail.

\subsection{Auxiliary basis functions}

Now, we present the construction of the auxiliary basis functions.
For each coarse element $K_{i}$, we define $V(K_i) = H^1(K_i)$, and
solve the spectral problem: find $(\lambda^{(i)}_j, \phi_{j}^{(i)}) \in \mathbb{R}\times V(K_{i})$
\[
a_i(\phi_{j}^{(i)},v)=\lambda_{j}^{(i)}s_i(\phi_{j}^{(i)},v), \quad\;\forall v\in V(K_{i}),
\]
where $a_i(u,v)=\int_{K_i}\kappa\nabla u\cdot\nabla v$, $s_i(u,v)=\int_{K_i} \tilde{\kappa} uv$, $\tilde{\kappa} = \kappa\sum_{j=1}^{N_c} |\nabla\chi_{j}|^{2}$
and  $\{\chi_j\}$ is a set of partition of unity functions with respect to the coarse grid.
We remark that one can take $\{ \chi_j\}$ as the
standard multiscale basis functions or the standard piecewise linear functions.
We assume that the eigenfunctions satisfy the normalized condition $s_i(\phi_j^{(i)},\phi_j^{(i)})=1$.
We can assume that the eigenvalues are sorted in ascending order, that is, $\lambda_{1}\leq\lambda_{2}\leq\dots$. We then choose the first $J_i$ eigenfunctions with small eigenvalue, and define the local auxiliary space $V_{aux}(K_{i})$ as the span of these eigen-basis functions, which is
\[
V_{aux}(K_{i})=\text{span}\{\phi_{j}^{(i)}|1\leq j\leq J_i\}.
\]
Notice that, by construction, $\lambda_{J_i+1}^{(i)} = O(1)$. We define
\begin{equation*}
\Lambda = \min_{1\leq i\leq N} \lambda_{J_i+1}^{(i)}.
\end{equation*}
Finally, the auxiliary space $V_{aux}$ is defined as the sum of these local auxiliary spaces, that is, $V_{aux}=\oplus_{i}V_{aux}(K_{i})$.

Using the auxiliary space, we define a projection operator $\pi:V\rightarrow V_{aux}$ by
\[
\pi(v)=\sum_{i=1}^N \sum_{j=1}^{J_i} s_{i}(v,\phi_{j}^{(i)})   \phi_{j}^{(i)}, \quad \;\forall v\in V.
\]
We also denote the kernel of the operator $\pi$ as
\[
\tilde{V}=\{w\in V|\;\pi(w)=0\}.
\]

\subsection{Multiscale basis functions}

\begin{figure}[ht]
\centering
\includegraphics[scale=0.5]{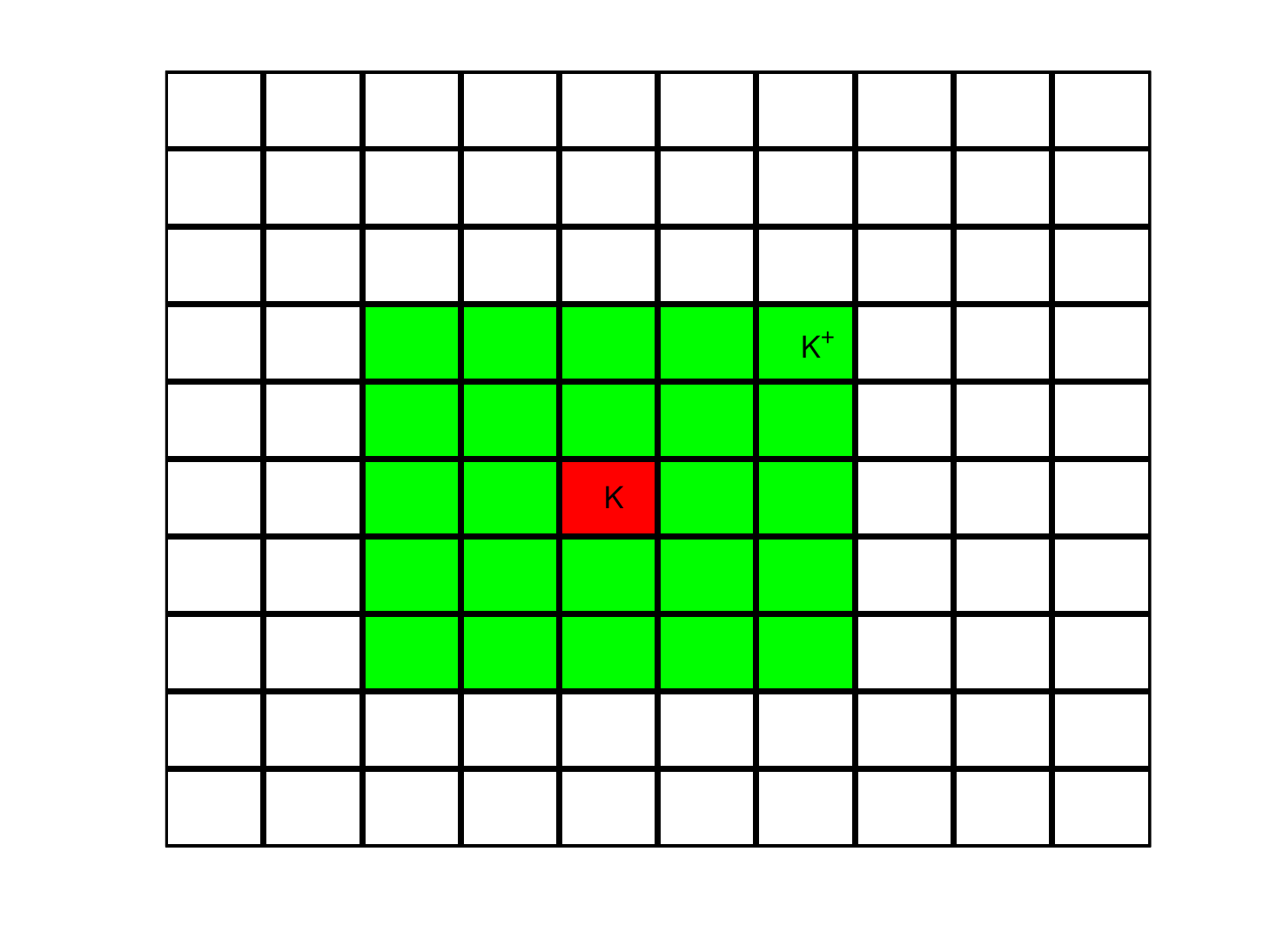}
\caption{A coarse element $K$ and its oversampled region $K^+$.}
\label{fig:kplus}
\end{figure}

Now we present the construction of offline multiscale basis functions using the auxiliary space.
For each coarse element $K_i$, we define an oversampled region $K_i^+$ by extending $K_i$ by $\ell$ coarse grid layers.
See Figure~\ref{fig:kplus} for an illustration of $K^+$ with $\ell = 2$.

For each $\phi_{j}^{(i)} \in V_{aux}$, we will construct a basis function $\psi^{(i)}_{j,ms}$ whose support is $K_i^+$.
Using the results in \cite{chung2017constraint},
the multiscale basis function $\psi^{(i)}_{j,ms}$ is constructed by solving the following local problem: find $\psi^{(i)}_{j,ms}\in V_0(K^+_i)$ such that
\begin{equation}
\label{eq:local_basis}
a(\psi_{j,ms}^{(i)},v)+s(\pi(\psi_{j,ms}^{(i)}),\pi(v))=s(\phi_{j}^{(i)},\pi(v)), \quad \;\forall v\in V_0(K^+_i),
\end{equation}
where $V_0(K_i^+) = H^1_0(K_i^+)$ and $s(u,v) = \sum_{i=1}^N s_i(u,v)$.
We remark that the above definition is defined in the continuous space $V_0(K_i^+)$.
In our numerical simulations, we solve the above problem using the fine mesh defined in $K_i^+$
and an appropriate finite element method.
Finally,
the multiscale finite element space $V_{ms}$ is defined as the span of these multiscale basis functions, namely, $V_{ms}=\text{span}\{\phi^{(i)}_{j,ms}\}$.
This method is called the CEM-GMsFEM.

We remark that our multiscale basis functions are used to approximate the related global basis functions.
The global basis function $\psi^{(i)}_{j}$ is defined by solving the following problem: find $\psi^{(i)}_{j}\in V$ such that
\begin{equation}
\label{eq:global_basis}
a(\psi_{j}^{(i)},v)+s(\pi(\psi_{j}^{(i)}),\pi(v))=s(\phi_{j}^{(i)},\pi(v)), \quad \;\forall v\in V.
\end{equation}
The global space is defined by $V_{glo}=\text{span}\{\psi_{j}^{(i)}\}$.
We note that these global basis functions have an exponential decay property (see \cite{chung2017constraint}),
which motivates the definitions of the multiscale basis functions $\psi^{(i)}_{j,ms}$ having local supports.
One important property of the global space is the orthogonal decomposition $V=V_{glo} \oplus \tilde{V}$
with respect to the inner product $a(u,v)$.
This global space will be used in the analysis of the convergence result of our online adaptive enrichment method.


\section{Online basis functions and adaptive enrichment}
\label{sec:online}
In this section, we will introduce an online enrichment method for this CEM. We will first present the construction of online basis functions. Then, we will give an error estimate for using the online enrichment method.
Here, for online basis function, we mean the basis functions constructed in online stage by using the residual of the solution which contain the information of the source term. We will construct the online basis function in an iterative process. We remark that the error will decay rapidly such that the error will within a acceptable range in the first or two iterations.

To begin, we define a residual functional $r:V \rightarrow \mathbb{R}$. Let $u_{ms}\in V_{ms}$ be a numerical solution computed by solving (\ref{eq:global_ms}). The residual functional $r$ is defined by
\[
r(v) = a(u_{ms},v) - \int_{\Omega} f v, \quad \forall v\in V.
\]
For each coarse neighborhood $\omega_i$, we define the local residual functional $r_i: V\rightarrow \mathbb{R}$ by
\[
r_i(v) = r(\chi_i v), \quad\forall v\in V.
\]
The local residual $r_i$ gives a measure of the error $u-u_{ms}$ in the coarse neighborhood $\omega_i$.

The construction of online basis function is related to the local residual $r_i$.
Using the local residual $r_i$, we will construct an online basis function $\phi_{on}^{(i)}$ whose support is an oversampled region $\omega_i^+$.
More precisely, the online basis function $\beta_{ms}^{(i)}\in V_0(\omega^+_i)$ is obtained by solving
\begin{equation}
\label{eq:on_basis}
a(\beta_{ms}^{(i)},v) + s(\pi(\beta_{ms}^{(i)}),\pi (v)) = r_i(v), \;\;\forall v\in V_0(\omega^+_i),
\end{equation}
where $V_0(\omega^+_i) = H^1_0(\omega_i^+)$.
We can perform the above construction for each $r_i$, or for some selected $r_i$ (with $i\in I$ for an index set $I$) based on an adaptive criterion.
We remark that the above online basis is obtained in the local region $\omega_i^+$. This is the result of a localization result
for the corresponding global online basis function $\beta_{glo}^{(i)} \in V$ defined by
\begin{equation}
\label{eq:global_on_basis}
a(\beta_{glo}^{(i)},v) + s(\pi(\beta_{glo}^{(i)}),\pi (v)) = r_i(v), \;\;\forall v\in V.
\end{equation}

After constructing the online basis functions, we can enrich our multiscale space by adding these online basis to the multiscale space, namely, $V_{ms}=V_{ms} + \text{span}_{i\in I}\{\beta_{ms}^{(i)}\}$. Using this multiscale finite element space, we can compute a new numerical solution by solving the equation (\ref{eq:global_ms}).
We can repeat the process to enrich our multiscale space until the residual norm is smaller than a given tolerance.
Next, we present the precise online adaptive enrichment algorithm.

\vspace{0.5cm}
\noindent
{\bf Online adaptive enrichment algorithm}

We first choose an initial space $V_{ms}^{(1)}$. This is the space obtained by using the offline basis functions
constructed in Section \ref{sec:method}.
We also choose a real number $\theta$ such that $0 \leq \theta < 1$. This number determines how many online basis functions
are needed in each online iteration.
Then, we will generate a sequence of spaces $V_{ms}^{(m)}$
and a sequence of multiscale solutions $u_{\text{ms}}^{(m)}$ obtained by solving
\eqref{eq:global_ms}.

For each $m= 1, 2, \dots$, we assume that $V_{ms}^{(m)}$ is given.
We will preform the following procedures to obtain the new multiscale space $V_{ms}^{(m+1)}$.

\begin{enumerate}
\item[Step 1:] Find the multiscale solution in the space $V^{(m)}_{ms}$. That is,
find $u_{ms}^{(m)} \in V^{(m)}_{ms}$ such that
\begin{equation}
a(u^{(m)}_{ms}, v) = (f, v), \quad \text{for all} \,\,\, v \in V^{(m)}_{ms}.
\label{eq:solve_u}
\end{equation}

\item[Step 2:] Compute the local residuals $z_i(v)$ where
\begin{equation*}
z_i(v) = a(u^{(m)}_{ms}, v) - (f, v), \quad \forall v\in V_0(\omega_i).
\end{equation*}
Define $\delta_i = \|z_i\|_{a^*}$ where $\|z_i\|_{a^*}=\sup_{v\in V_0(\omega_i)} \cfrac{r(v)}{\|v\|_a}$.
We re-numerate the indices of $\omega_i$ such that
$\delta_1 \geq \delta_2 \geq \cdots$.
Choose the first $k$ regions so that
\begin{equation}
\label{eq:theta}
\sum_{i=k+1}^N \delta_i^2 < \theta \sum_{i=1}^N \delta_i^2.
\end{equation}

\item[Step 3:] Compute the local online basis functions. For each $1\leq i\leq k$ and
coarse region $\omega_i$,
 find $\beta^{(i)}_{ms}\in V_0(\omega_i^+)$ such that
\begin{equation*}
a(\beta_{ms}^{(i)},v) + s(\pi(\beta_{ms}^{(i)}),\pi (v)) = r^{(m)}_i(v) \;\;\forall v\in V_0(\omega^+_i)
\end{equation*}
where $r^{(m)}_i(v) = a(u^{(m)}_{ms},\chi_i v) - \int_{\Omega} f \chi_i v$.
\item[Step 4:] Enrich the multiscale space. Let
\[
V^{(m+1)}_{ms} = V^{(m)}_{ms} + \text{span}_{1\leq i \leq k} \{ \beta^{(i)}_{ms}\}.
\]
\end{enumerate}

In the next section, we will analyze the convergence rate for this online adaptive enrichment method.
In particular, we will prove the following theorem.

\begin{theorem}
\label{thm:conv}
Let $u$ be the solution of (\ref{eq:original}) and let $\{ u_{ms}^{(m)}\}$ be
the sequence of multiscale solutions obtained by our online adaptive enrichment algorithm. Then
we have
\[
\|u-u_{ms}^{(m+1)}\|^2_{a}\leq 3(1+\Lambda^{-1})\Big(C(\ell+1)^{d}E + 2M^2\theta \Big)\|u-u^{(m)}_{ms}\|^2_a
\]
where $E=3(1+\Lambda^{-1})(1+2(1+\Lambda^{\frac{1}{2}})^{-1})^{1-\ell}$,
$M$ is maximum number of overlapping subdomains
and $C$ is a constant.
\end{theorem}

%

\begin{remark}
We note that the convergence rate depends two terms $C(\ell+1)^{d}E$ and $2M^2\theta$.
By using enough number of oversampling layers, the term $C(\ell+1)^{d}E$ tends to zero.
Thus, the factor $2M^2\theta$ dominates the convergence rate. One can choose $\theta$
to obtain a desired convergence rate. We will also confirm this by some numerical examples.
This is an improvement over the online method in \cite{chung2015residual},
where the convergence rate is $(C_1 + C_2 \theta)$ with $0 < C_1 < 1$.
\end{remark}

\section{Numerical Result}
\label{sec:num}

\begin{figure}[ht!]
\centering
\includegraphics[scale=0.5]{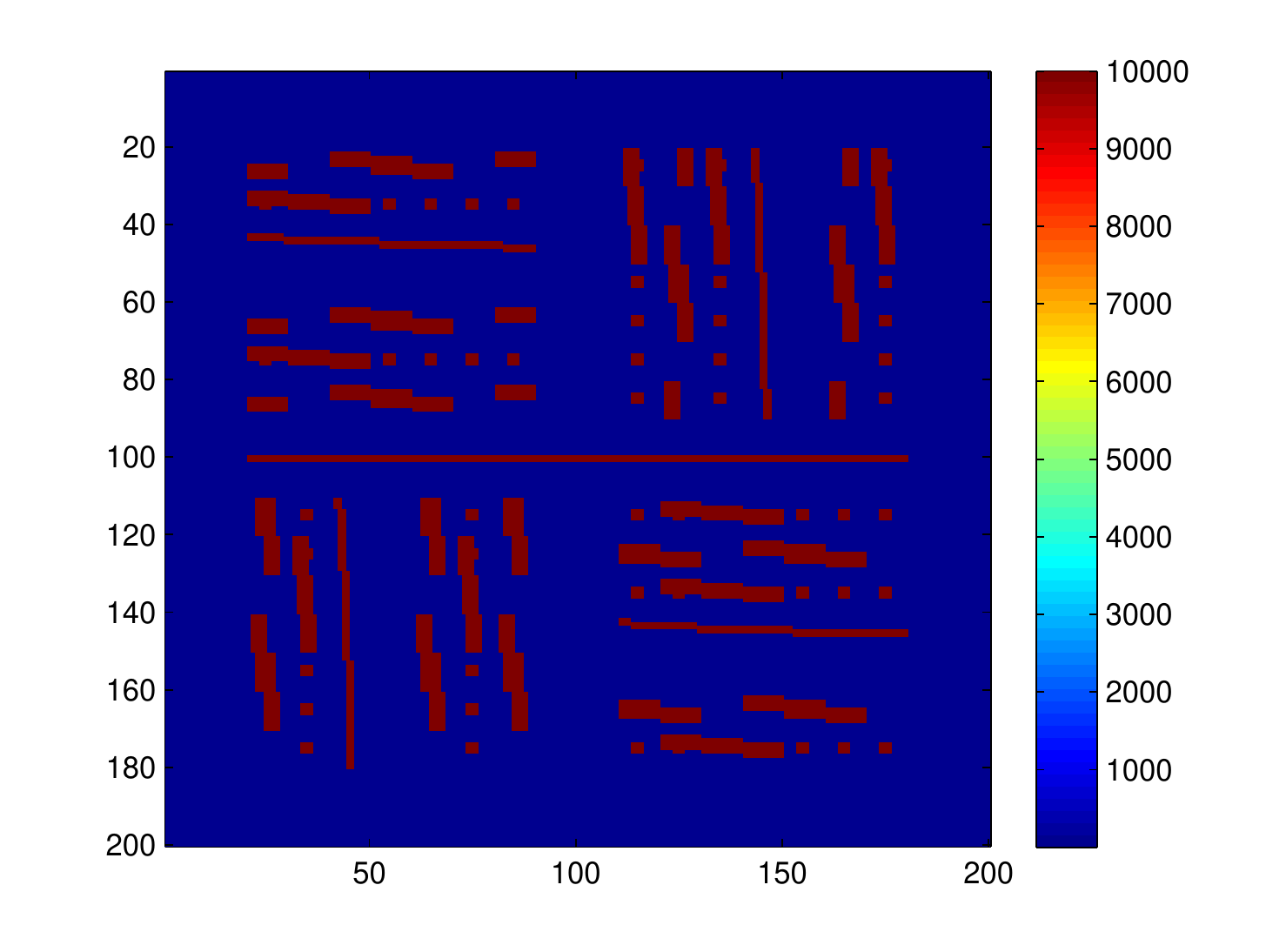}
\caption{The medium $\kappa$ for the test case $1$.}
\label{fig:case1_kappa}
\end{figure}

In this section, we present some numerical results to demonstrate the convergence of our proposed method.
We take the computational domain $\Omega = (0,1)^2$.
The medium parameter $\kappa$ in the equation (\ref{eq:original})
is chosen to be the function shown in Figure \ref{fig:case1_kappa}.
We note that the medium $\kappa$ contains high contrast inclusions and channels.
The fine mesh size $h$ is taken to be $1/200$, while the coarse mesh size $H$ in this example is $1/10$.
In all our results, we take the number of oversampling layers $\ell=2$.
We will illustrate the performance of our method by using two different source terms
$f_1=((x-0.5)^2+(y-0.5)^2)^{-\frac{1}{4}}$ and $f_2=((x-0.5)^2+(y-0.5)^2)^{-\frac{3}{4}}$.
We will test the performance by considering uniform enrichment
and by using the online adaptive enrichment algorithm presented in Section \ref{sec:online}.


In Table \ref{tab:case1_error1}, we present the $L_{2}$ error and the energy error
for the case $f_1$ with uniform enrichment, that is $\theta=0$.
From the first two online iterations, we observe very fast convergence
of the method.
Next, we will consider some adaptive results for this case.
In Table \ref{tab:case1_error1_ad_095}, we present the error decay
by using our online adaptive enrichment algorithm with $\theta=0.95$.
That is, we only
add basis for regions which account for the largest $5\%$ of the residual.
From the table, we observe that the convergence rate in the energy norm
is $0.9154$, which is close to $0.95$.
This results confirm our assertion that the convergence rate can be controlled by the user-defined parameter $\theta$.
We remark that the convergence rate is computed by taking the maximum of
all $\|u - u_{ms}^{(m+1)}\|_a^2 / \| u-u_{ms}^{(m)}\|_a^2$.
In Table \ref{tab:case1_error2_ad_01}, we present the adaptive result with $\theta=0.1$.
That is, we
add basis for regions which account for the largest $90\%$ of the residual.
From the table, we observe that the convergence rate in the energy norm
is $0.0589$, which is close to $0.1$.
This result also confirms our prediction.
Moreover, we note that the adaptive approach allows adding a very few
online basis functions to reduce the error to $1\%$.

\begin{table}[ht]
\centering
\begin{tabular}{|c|c|c|c|c|}
\hline
Number of offline basis & online iteration & oversampling layers & $L_{2}$ error & energy error\tabularnewline
\hline
3 & 0 & 2 & 0.37\% & 4.71\%\tabularnewline
\hline
3 & 1 & 2 & 6.75e-05\% & 1.28e-03\%\tabularnewline
\hline
3 & 2 & 2 & 1.57e-08\% & 2.64e-08\%\tabularnewline
\hline
\end{tabular}

\caption{Using source term $f_1$ and uniform enrichment.}
\label{tab:case1_error1}
\end{table}

\begin{table}[ht]
\centering
\begin{tabular}{|c|c|c|c|c|}
\hline
Number of offline basis $\omega$ & DOF & oversampling layers & $L_{2}$ error & energy error\tabularnewline
\hline
3 & 300 & 2 & 0.37\% & 4.71\%\tabularnewline
\hline
3 & 311 & 2 & 0.14\% & 2.21\%\tabularnewline
\hline
3 & 339 & 2 & 0.073\% & 1.12\%\tabularnewline
\hline
3 & 368 & 2 & 0.033\% & 0.57\%\tabularnewline
\hline
\end{tabular}
\caption{Using source term $f_1$ and online adaptivity with $\theta=0.95$. Convergence rate is $0.9154$.}
\label{tab:case1_error1_ad_095}
\end{table}

\begin{table}[ht]
\centering
\begin{tabular}{|c|c|c|c|c|}
\hline
Number of offline basis & DOF & oversampling layers & $L_{2}$ error & energy error\tabularnewline
\hline
3 & 300 & 2 & 0.37\% & 4.71\%\tabularnewline
\hline
3 & 341 & 2 & 0.073\% & 1.09\%\tabularnewline
\hline
3 & 407 & 2 & 0.014\% & 0.21\%\tabularnewline
\hline
3 & 470 & 2 & 2.93e-03\% & 0.051\%\tabularnewline
\hline
\end{tabular}
\caption{Using source term $f_1$ and online adaptivity with $\theta=0.1$. Convergence rate is $0.0589$.}
\label{tab:case1_error2_ad_01}
\end{table}

%

%
%

\begin{figure}[ht!]
\centering
\includegraphics[scale=0.5]{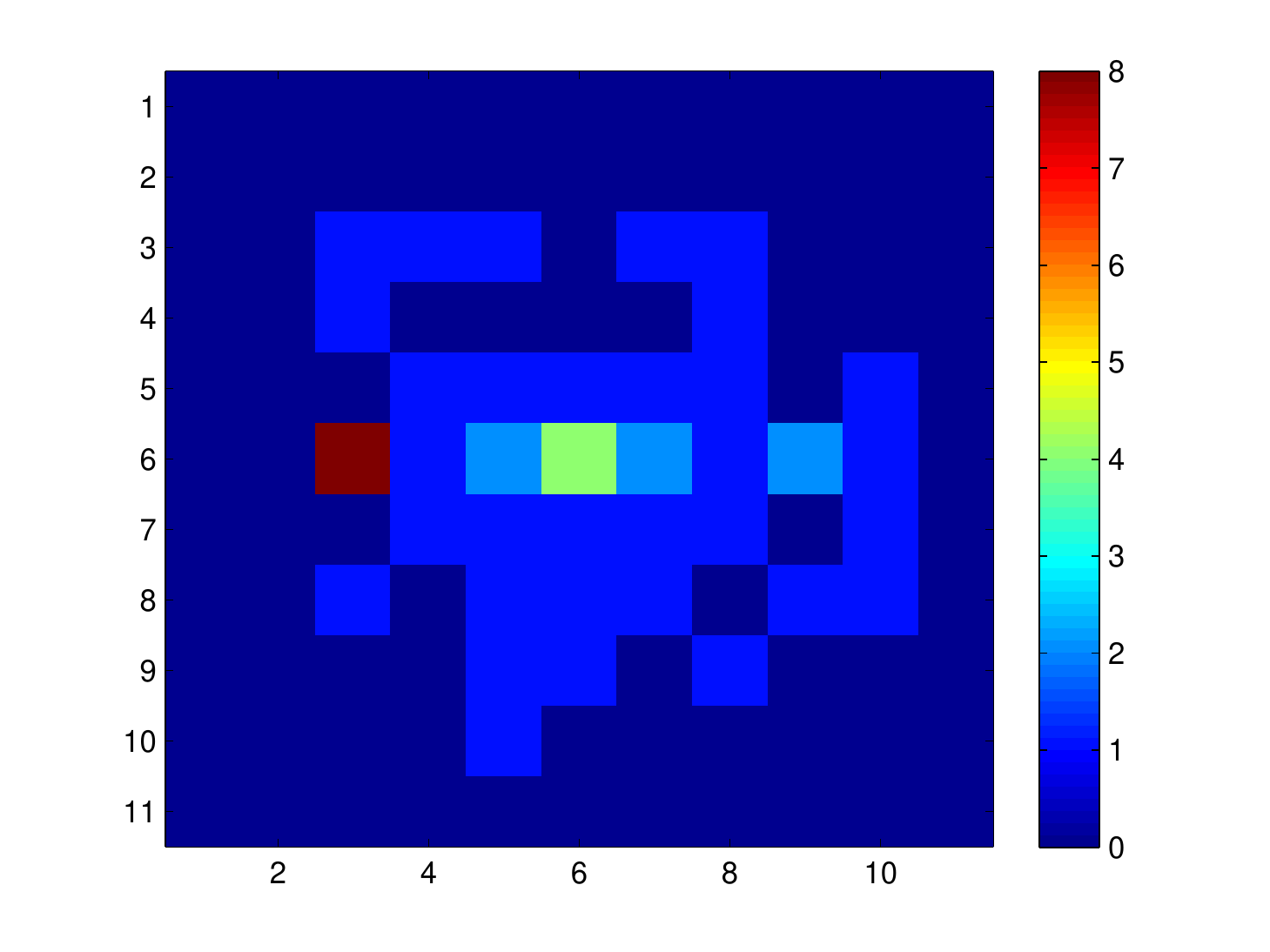}
\includegraphics[scale=0.5]{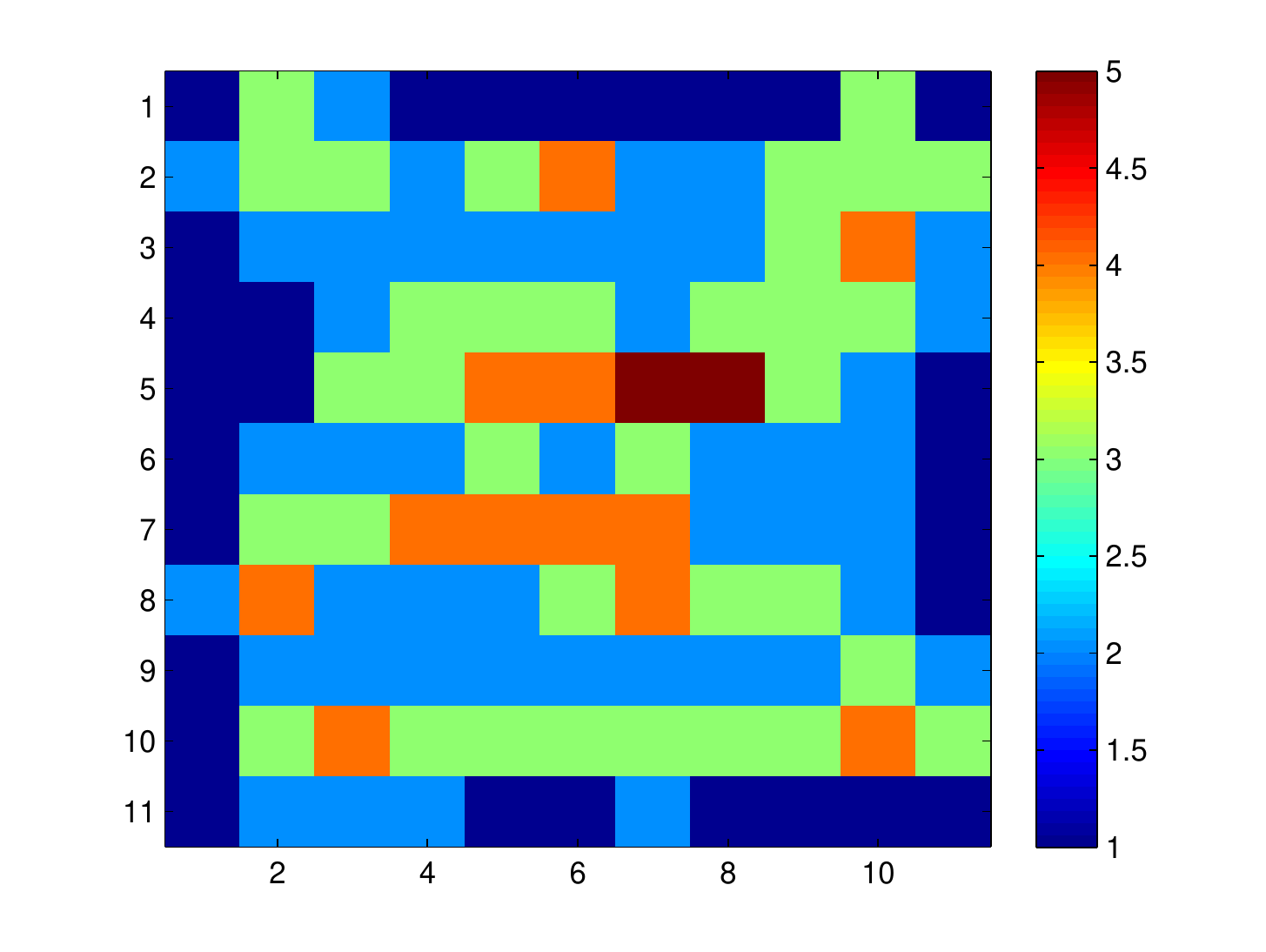}
\caption{Number of online basis functions for the source $f_2$: Left: $\theta=0.95$. Right: $\theta=0.1$.}
\label{fig:numbasis}
\end{figure}

Now, we consider the second source term $f_2$.
In Table \ref{tab:case2_error2}, we present the error decay using uniform enrichement.
We observe very fast decay in error from this table.
Next, we test the performance using adaptivity.
In Tables \ref{tab:case2_error2_ad_095} and \ref{tab:case2_error2_ad_01},
we present the error decays with $\theta=0.95$ and $\theta=0.1$ respectively.
We observe that the convergence rates in these two cases are $0.9338$ and $0.09$ respectively.
This confirms that the user-defined parameter is useful in controlling the convergence rate
of our adaptive method.
Moreover, we note that the adaptive approach allows adding a very few
online basis functions to reduce the error to $1\%$.
Furthermore, in Figure \ref{fig:numbasis}, we show the number of online basis functions
added in the computational domain.
For $\theta=0.95$, we will add a small number of basis functions in each iteration. We observe that
the basis functions are added near the singularity of the source $f_2$ and along the high contrast channel.
For $\theta=0.1$, more basis functions are added throughout the domain with a faster convergence rate.
We still observe that more basis are added near the singularity of $f_2$ and along the high contrast channel in $\kappa$.

\begin{table}[ht]
\centering
\begin{tabular}{|c|c|c|c|c|}
\hline
Number of offline basis & online iteration & oversampling layers & $L_{2}$ error & energy error\tabularnewline
\hline
3 & 0 & 2 & 1.06\% & 11.70\%\tabularnewline
\hline
3 & 1 & 2 & 6.43e-05\% & 1.51e-03\%\tabularnewline
\hline
3 & 2 & 2 & 1.57e-08\% & 4.25e-08\%\tabularnewline
\hline
\end{tabular}
\caption{Using source term $f_2$ and uniform enrichment.}
\label{tab:case2_error2}
\end{table}


\begin{table}[ht]
\centering
\begin{tabular}{|c|c|c|c|c|}
\hline
Number of offline basis & DOF & oversampling layers & $L_{2}$ error & energy error\tabularnewline
\hline
3 & 300 & 2 & 1.06\% & 11.70\%\tabularnewline
\hline
3 & 309 & 2 & 0.13\% & 1.95\%\tabularnewline
\hline
3 & 324 & 2 & 0.062\% & 0.99\%\tabularnewline
\hline
3 & 347 & 2 & 0.031\% & 0.51\%\tabularnewline
\hline
\end{tabular}
\caption{Using source term $f_2$ and online adaptivity with $\theta=0.95$. Convergence rate is $0.9338$.}
\label{tab:case2_error2_ad_095}
\end{table}

\begin{table}[ht]
\centering
\begin{tabular}{|c|c|c|c|c|}
\hline
Number of offline basis & DOF & oversampling layers  & $L_{2}$ error & energy error\tabularnewline
\hline
3 & 300 & 2 & 1.06\% & 11.70\%\tabularnewline
\hline
3 & 391 & 2 & 9.27e-03\% & 0.13\%\tabularnewline
\hline
3 & 513 & 2 & 3.65e-04\% & 6.44e-03\%\tabularnewline
\hline
3 & 578 & 2 & 5.31e-05\% & 9.89e-04\%\tabularnewline
\hline
\end{tabular}
\caption{Using source term $f_2$ and online adaptivity with $\theta=0.1$. Convergence rate is $0.09$.}
\label{tab:case2_error2_ad_01}
\end{table}

\begin{figure}[ht!]
\centering
\includegraphics[scale=0.5]{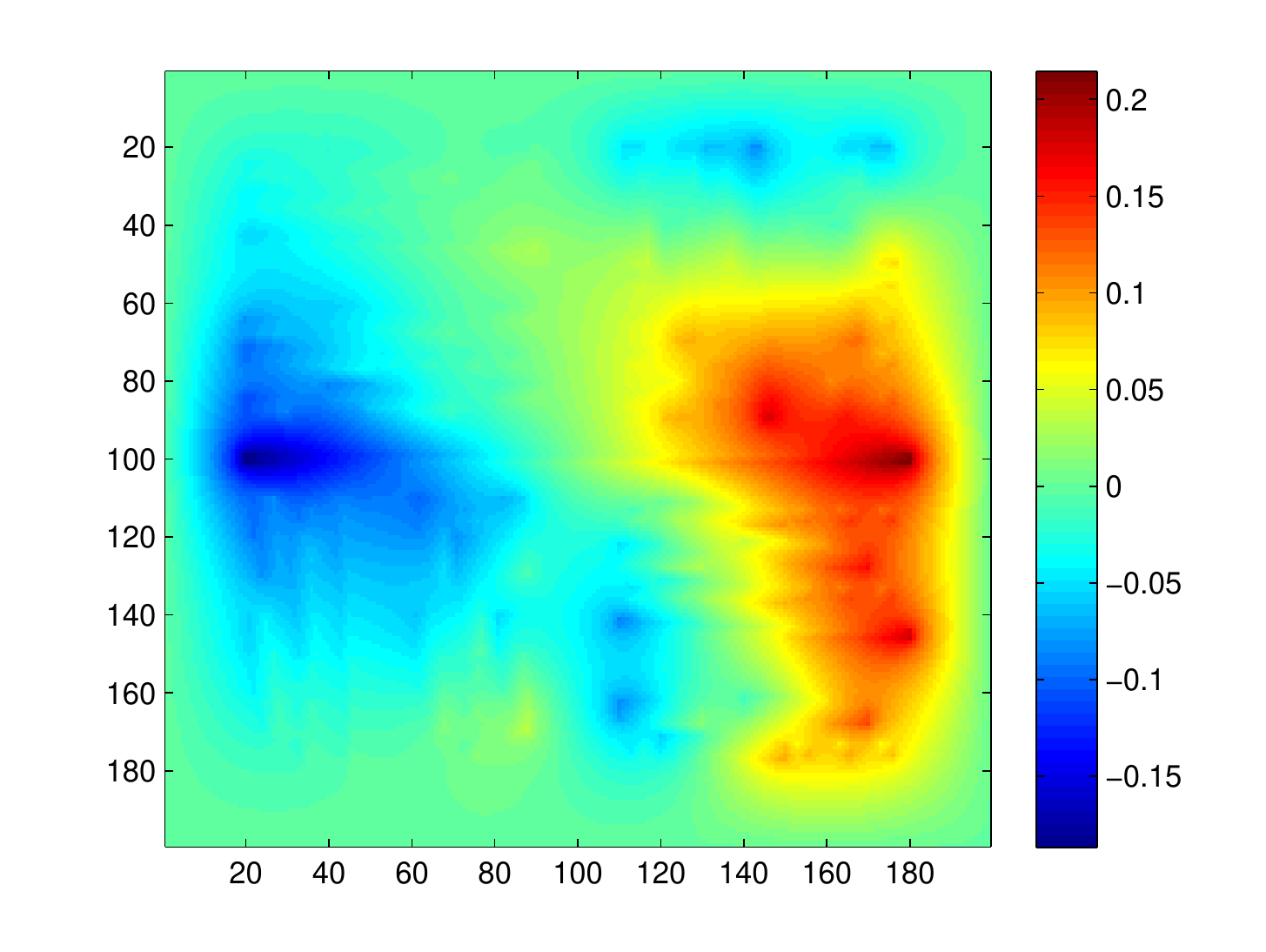}
\includegraphics[scale=0.5]{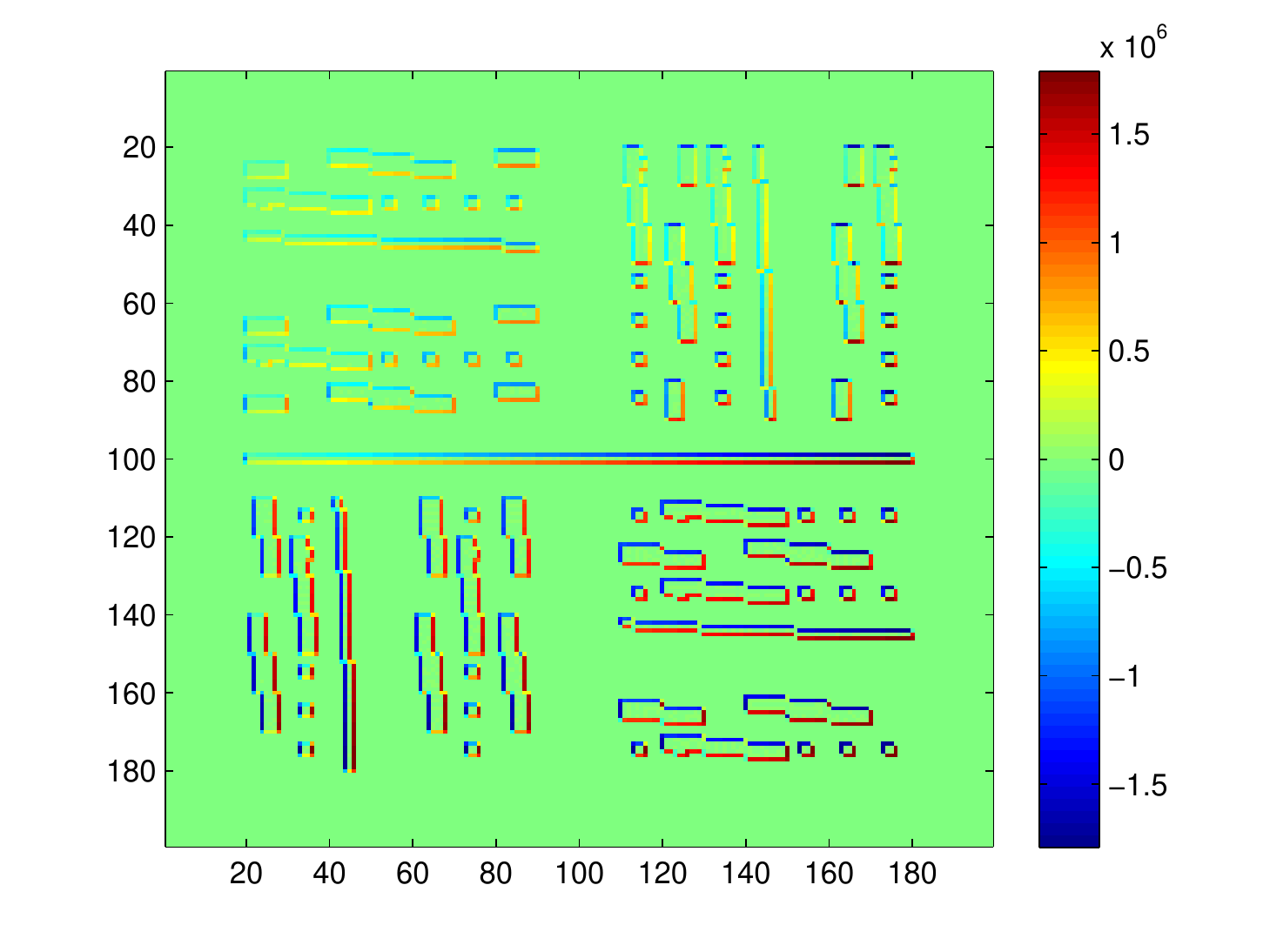}
\caption{Left: Reference solution for the source $f_3$. Right: the source term $f_3$.}
\label{fig:sol_case4}
\end{figure}

Finally, we present a test case with a more singular source term $f_3=-\nabla \cdot(\kappa \nabla (xy))$,
shown in Figure \ref{fig:sol_case4}, where the reference solution is also presented.
In Table \ref{tab:case4_error1}, we present the error decay with uniform enrichment
and observe the same type of exponential decay as the earlier examples.
We also observe that the error is relatively large when no online basis function is used.
In Table \ref{tab:case4_error2_ad_01}, we present the results with the online adaptive enrichment algorithm with $\theta=0.1$.
We see that the numerically computed convergence rate is $0.0771$, which is close to the parameter $\theta$.

\begin{table}[ht]
\centering
\begin{tabular}{|c|c|c|c|c|}
\hline
Number of offline basis & online iteration & oversampling layers & $L_{2}$ error & energy error\tabularnewline
\hline
3 & 0 & 2 & 30.01\% & 82.57\%\tabularnewline
\hline
3 & 1 & 2 & 0.0066\% & 0.0030\%\tabularnewline
\hline
3 & 2 & 2 & 4.45e-07\% & 1.22e-07\%\tabularnewline
\hline
\end{tabular}

\caption{Using source term $f_3$ and uniform enrichment.}
\label{tab:case4_error1}
\end{table}

\begin{table}[ht]
\centering
\begin{tabular}{|c|c|c|c|c|}
\hline
Number of offline basis & DOF & oversampling layers  & $L_{2}$ error & energy error\tabularnewline
\hline
3 & 300 & 2 & 30.01\% & 82.57\%\tabularnewline
\hline
3 & 356 & 2 & 8.68\% & 22.06\%\tabularnewline
\hline
3 & 378 & 2 & 4.87\% & 5.41\%\tabularnewline
\hline
3 & 392 & 2 & 4.46\% & 1.50\%\tabularnewline
\hline
\end{tabular}
\caption{Using source term $f_3$ and online adaptivity with $\theta=0.1$. Convergence rate is $0.0771$.}
\label{tab:case4_error2_ad_01}
\end{table}

\section{Convergence analysis}
\label{sec:analysis}

In this section, we analyze the convergence of the online adaptive enrichment algorithm presented in Section \ref{sec:online}.
First, we need some notations.
We will define two different norms for the space $V$.
One is the $a$-norm $\|\cdot\|_a$ where $\|u\|_a^2=\int_\Omega \kappa | \nabla u|^2$. The other is $s$-norm $\|\cdot\|_s$ where $\|u\|_s^2=\int_\Omega \tilde{\kappa} u^2$.
For a given subdomain $\Omega_i \subset \Omega$, we will define the local $a$-norm and $s$-norm by $\|u\|_{a(\Omega_i)}^2=\int_{\Omega_i} \kappa | \nabla u|^2$ and $\|u\|_{s(\Omega_i)}^2=\int_{\Omega_i} \tilde{\kappa}  u^2$ respectively.

Next, we will recall a few theoretical results from \cite{chung2017constraint}
that are useful for our analysis.
The first result is Lemma \ref{lem:infsup}.
\begin{lemma}
\label{lem:infsup}
There is a constant $D$ such that
for all $v_{aux}\in V_{aux}$
there exists a function $v\in V$ such that
\[
\pi(v)=v_{aux},\qquad\|v\|_{a}^{2}\leq D \, \|v_{aux}\|_{s}^{2},\qquad\text{supp}(v)\subset\text{supp}(v_{aux}).
\]
\label{lem:projection}
\end{lemma}

The second result is a localization result, saying that the global basis function
defined in (\ref{eq:global_basis}) has an exponential decay
outside an oversampled region. This result motivates the local multiscale basis functions defined in (\ref{eq:local_basis}).

\begin{lemma}
\label{lem:new}
We consider the oversampled domain $K_i^+$ obtained from $K_i$ by extending $\ell$ coarse grid layers
with $\ell\geq 2$.
Let $\phi_j^{(i)} \in V_{aux}$ be a given auxiliary multiscale basis function.
We let $\psi_{j,ms}^{(i)}$ be the multiscale basis functions obtained in (\ref{eq:local_basis})
and let $\psi_{j}^{(i)}$ be the global multiscale basis functions obtained in (\ref{eq:global_basis}).
Then we have
\[
\|\psi_{j}^{(i)}-\psi_{j,ms}^{(i)}\|_{a}^{2}+\|\pi(\psi_{j}^{(i)}-\psi_{j,ms}^{(i)})\|_{s}^{2}\leq E
\Big(\|\psi_{j}^{(i)}\|_{a}^{2}+\|\pi(\psi_{j}^{(i)})\|_{s}^{2}\Big)
\]
where $E = 3(1+\Lambda^{-1}) \Big (1+(2(1+\Lambda^{-\frac{1}{2}}))^{-1}\Big)^{1-\ell}$.
\end{lemma}

Next, we will need the following lemma in our analysis. The proof is given in the Appendix.
\begin{lemma}
\label{lem:sum}
Assume the same conditions in Lemma \ref{lem:new}. We have
\begin{align*}
&\; \| \sum_{i=1}^N \sum_{j=1}^{J_i} c^{(i)}_j( \psi_j^{(i)} - \psi_{j,ms}^{(i)} )\|_a^2 + \| \pi(\sum_{i=1}^N \sum_{j=1}^{J_i} c^{(i)}_j( \psi_j^{(i)} - \psi_{j,ms}^{(i)} ))\|_s^2\\
\leq & \; C (1+\Lambda^{-1}) (\ell+1)^{d} \sum_{i=1}^N  \Big(\| \sum_{j=1}^{J_i}c^{(i)}_j(\psi_j^{(i)} - \psi_{j,ms}^{(i)})\|_a^2+\| \pi(\sum_{j=1}^{J_i}c^{(i)}_j(\psi_j^{(i)} - \psi_{j,ms}^{(i)}))\|_s^2\Big).
\end{align*}
\end{lemma}

In the next lemma, we see that the same localization result
holds for the online basis function defined in (\ref{eq:on_basis}). The
proof of the following lemma is the same as that for Lemma \ref{lem:new} and Lemma \ref{lem:sum},
and is therefore omitted.

\begin{lemma}
\label{lem:on}
We consider the oversampled domain $\omega_i^+$ obtained from $\omega_i$ by extending $\ell$ coarse grid layers
with $\ell\geq 2$.
We let $\beta_{ms}^{(i)}$ be the online multiscale basis functions obtained in (\ref{eq:on_basis})
and let $\beta_{glo}^{(i)}$ be the global online multiscale basis functions obtained in (\ref{eq:global_on_basis}).
Then we have
\[
\|\beta_{glo}^{(i)}-\beta_{ms}^{(i)}\|_{a}^{2}+\|\pi(\beta_{glo}^{(i)}-\beta_{ms}^{(i)})\|_{s}^{2}\leq E \Big(\|\beta_{glo}^{(i)}\|_{a}^{2}+\|\pi(\beta_{glo}^{(i)})\|_{s}^{2}\Big)
\]
where $E = 3(1+\Lambda^{-1}) \Big (1+(2(1+\Lambda^{-\frac{1}{2}}))^{-1}\Big)^{1-\ell}$. Furthermore, we have
\begin{align*}
&\; \| \sum_{i=1}^{N_c} ( \beta_{glo}^{(i)} - \beta_{ms}^{(i)} )\|_a^2 + \| \pi(\sum_{i=1}^{N_c} ( \beta_{glo}^{(i)} - \beta_{ms}^{(i)} ))\|_s^2\\
\leq & \; C (1+\Lambda^{-1}) (\ell+1)^{d} \sum_{i=1}^{N_c}  \Big(\| \beta_{glo}^{(i)} - \beta_{ms}^{(i)}\|_a^2+\| \pi(\beta_{glo}^{(i)} - \beta_{ms}^{(i)})\|_s^2\Big).
\end{align*}
\end{lemma}

Finally, we define a constant $C_0 \in \mathbb{R}$ as
\[
C_0 = \sup_{v \in V}\cfrac{\| \pi(v) \|^2_s}{\|v\|^2_a}.
\]
We remark that $C_0\leq \max\{\tilde{\kappa}\}C^2_{p}$, where $C_p$ is the Poincare constant defined
by
$\|w\|_{L^2(\Omega)} \leq C_p \| \nabla w\|_{L^2(\Omega)}$ for $w\in H^1_0(\Omega)$.


Now we are ready to prove Theorem \ref{thm:conv}.

\subsection{Proof of Theorem \ref{thm:conv}}

First, by using the Galerkin orthogonality, we have
\begin{equation}
\label{eq:galerkin}
\| u - u_{ms}^{(m+1)}\|_a \leq \| u - w\|_a, \quad \forall w\in V_{ms}^{(m+1)}.
\end{equation}
The proof is based on a suitable choice of $w\in V_{ms}^{(m+1)}$,
and consists of $4$ steps.

\noindent
{\bf Step 1}:

In this step, we will give a representation of the error $u-u_{ms}^{(m)}$.
For each $i=1,2,\cdots, N_c$,
we construct a global online basis function $\beta_{glo}^{(i)}\in V$ such that
\[
a(\beta_{glo}^{(i)},v)+s(\pi(\beta_{glo}^{(i)}),\pi(v))=r_i^{(m)}(v), \quad \;\forall v\in V.
\]
Summing over all $i=1,2,\cdots, N_c$, we have
\begin{align*}
a(\sum_{i=1}^{N_c}\beta_{glo}^{(i)},v)+s(\pi(\sum_{i=1}^{N_c}\beta_{glo}^{(i)}),\pi(v)) & =\sum_{i=1}^{N_c} r^{(m)}_i(v).
\end{align*}
By the definition of $r_i^{(m)}(v)$, we have
\begin{equation*}
\sum_{i=1}^{N_c} r^{(m)}_i(v) = \sum_{i=1}^{N_c} \Big( a(u_{ms}^{(m)}, \chi_i v) -\int_{\Omega} f\chi_i v \Big)
= \sum_{i=1}^{N_c} \Big( a(u_{ms}^{(m)}, \chi_i v) - a(u, \chi_i v) \Big)
=  a(u_{ms}^{(m)}-u,v).
\end{equation*}
Therefore, we have
\begin{equation}
\label{eq:rep2}
a(u-u_{ms}^{(m)}+\sum_{i=1}^{N_c}\beta_{glo}^{(i)},v)=s(-\pi(\sum_{i=1}^{N_c}\beta_{glo}^{(i)}),\pi(v)), \quad \forall v\in V.
\end{equation}
From the above relation, we see that
\begin{equation*}
a(u-u_{ms}^{(m)}+\sum_{i=1}^{N_c}\beta_{glo}^{(i)},v)= 0, \quad \forall v\in \tilde{V}.
\end{equation*}
Using the decomposition $V=V_{glo} \oplus \tilde{V}$, we have
\[
u-u_{ms}^{(m)}+\sum_{i=1}^{N_c}\beta_{glo}^{(i)}\in\tilde{V}^{\perp}=V_{glo}.
\]
Hence, we obtain the representation
\begin{equation}
\label{eq:rep}
u-u_{ms}^{(m)}+\sum_{i=1}^{N_c}\beta_{glo}^{(i)}=\sum_{i=1}^N\sum_{j=1}^{J_i} c_{j}^{(i)}\psi_{j}^{(i)}
\end{equation}
where $c_j^{(i)}$ are some coefficients.
We will use this representation in the next steps.
We remark that, in Step 2 and Step 3,
we will localize the terms $\psi_j^{(i)}$ and $\beta_{glo}^{(i)}$,
and estimate the errors.


\noindent
{\bf Step 2}:

In this step, we will localize each $\psi_j^{(i)}$ in (\ref{eq:rep}) and estimate the error.
In particular,
we will estimate $\|\sum_{i=1}^N\sum_{j=1}^{J_i} c_{j}^{(i)}\psi_{j}^{(i)}-\sum_{i=1}^N \sum_{j=1}^{J_i}c_{j}^{(i)}\psi_{j,ms}^{(i)}\|_{a}$.
We define $\eta := u-u_{ms}^{(m)}+\sum_{i=1}^{N_c}\beta_{glo}^{(i)}$. Using (\ref{eq:rep}), we have
\begin{align}
a(\eta,v) + s(\pi (\eta),\pi(v)) &= a(\sum_{i=1}^{N_c}\sum_{j=1}^{J_i}c_{j}^{(i)}\psi_{j}^{(i)},v) + s(\sum_{i=1}^{N_c}\sum_{j=1}^{J_i}c_{j}^{(i)}\pi(\psi_{j}^{(i)}),\pi(v))\\
&= \sum_{i=1}^{N_c}\sum_{j=1}^{J_i} c_j^{(i)} s(\phi_{j}^{(i)},\pi(v)), \quad \forall v\in V, \label{eq:rep1}
\end{align}
where the last equality follows from (\ref{eq:global_basis}).
We let $v_{aux}^{(i)} = \sum_{j=1}^{J_i} c_{j}^{(i)} \phi_{j}^{(i)} \in V_{aux}(K_i)$.
By Lemma \ref{lem:infsup},
there exists a function $q^{(i)}\in V_0(K_i)$ such that $\pi(q^{(i)}) = v_{aux}^{(i)}$ and
\[
\|q^{(i)}\|^2_a \leq D \| v_{aux}^{(i)} \|^2_s.
\]
Taking $v=q^{(i)}$ in (\ref{eq:rep1}),
we have
\begin{align*}
\|v_{aux}^{(i)}\|^2_s  & = a(\eta,q^{(i)}) + s(\pi(\eta),\pi(q^{(i)}))\\
&\leq  \Big(\|\eta\|^2_{a(K_i)} + \|\pi(\eta)\|^2_{s(K_i)}\Big)^{\frac{1}{2}}\Big(\|q^{(i)}\|^2_a + \|\pi(q^{(i)})\|^2_s\Big)^{\frac{1}{2}}\\
& \leq \Big(\|\eta\|^2_{a(K_i)} + \|\pi(\eta)\|^2_{s(K_i)}\Big)^{\frac{1}{2}}\Big((1+D)\|v_{aux}^{(i)}\|^2_s\Big)^{\frac{1}{2}}.
\end{align*}
Thus, by the orthogonality of the eigenfunctions $\phi_j^{(i)}$ and the normalization condition $s_i(\phi_j^{(i)},\phi_j^{(i)})=1$, we have
\begin{align*}
\sum_{i=1}^N\sum_{j=1}^{J_i}(c_{j}^{(i)})^{2} = \sum_{i=1}^N \|\sum_{j=1}^{J_i}c_{j}^{(i)} \phi_{j}^{(i)}\|^2_s  \leq (1+D) \Big(\|\eta\|^2_{a} + \|\pi(\eta)\|^2_{s}\Big)
\leq (1+D) (C_0+1)\|\eta\|^2_{a}
\end{align*}
where the last inequality follows from the definition of the constant $C_0$.
Recalling the definition of $\eta$, we have
\begin{align*}
\sum_{i=1}^N \sum_{j=1}^{J_i}(c_{j}^{(i)})^{2} & \leq (1+D) (C_0+1)  \|u-u_{ms}^{(m)}+\sum_{i=1}^{N_c}\beta_{glo}^{(i)}\|_{a}^{2}\\
& \leq 2(1+D) (C_0+1)\Big(\|u-u_{ms}^{(m)}\|_{a}^{2}+\|\sum_{i=1}^{N_c}\beta_{glo}^{(i)}\|_{a}^{2}\Big)\\
 & \leq 4(1+D) (C_0+1)\|u-u_{ms}^{(m)}\|_{a}^{2}
\end{align*}
where the last inequality follows from (\ref{eq:rep2}).
Finally, using Lemma \ref{lem:new} and Lemma \ref{lem:sum}, we have
\[
\|\sum_{i=1}^N \sum_{j=1}^{J_i}c_{j}^{(i)}\psi_{j}^{(i)} - \sum_{i=1}^N \sum_{j=1}^{J_i} c_{j}^{(i)}\psi_{j,ms}^{(i)}\|_{a}^2
\leq C(\ell+1)^{d}E(1+D) (C_0+1)\|u-u_{ms}^{(m)}\|_{a}^{2}.
\]

\noindent
{\bf Step 3}:

In this step,
we will derive an estimate for $\beta_{glo}^{(i)}-\beta_{ms}^{(i)}$. By Lemma \ref{lem:on},
we have
\begin{align*}
\|\beta_{glo}^{(i)}-\beta_{ms}^{(i)}\|_{a}^{2} +\|\pi(\beta_{glo}^{(i)}-\beta_{ms}^{(i)})\|_{s}^{2} & \leq E\Big(\|\beta_{glo}^{(i)}\|_{a}^{2} + \|\pi(\beta_{glo}^{(i)})\|_{s}^{2}\Big).
\end{align*}
Using the equation (\ref{eq:global_on_basis}),
\begin{align*}
\|\beta_{glo}^{(i)}\|_{a}^{2} + \|\pi(\beta_{glo}^{(i)})\|_{s}^{2} &=  a(u-u^{(m)}_{ms},\chi_i \beta_{glo}^{(i)})\\
& \leq \|u-u^{(m)}_{ms}\|_{a(\omega_i)} \|\chi_i \beta_{glo}^{(i)}\|_{a(\omega_i)} \\
& \leq \sqrt{2}\|u-u^{(m)}_{ms}\|_{a(\omega_i)} \Big(\|\beta_{glo}^{(i)}\|^2_{a(\omega_i)} + \|\beta_{glo}^{(i)}\|^2_{s(\omega_i)}\Big)^{\frac{1}{2}}.
\end{align*}
Since $\pi$ is an orthogonal projection onto the space spanned by the eigenfunctions $\{ \phi_j^{(i)}\}$ for $j=1,2,\cdots, J_i$ and $i=1,2,\cdots, N$, we have
\begin{align*}
\|\beta_{glo}^{(i)}\|^2_{s(\omega_i)} &\leq \|\pi(\beta_{glo}^{(i)})\|^2_{s(\omega_i)} + \|(I-\pi)(\beta_{glo}^{(i)})\|^2_{s(\omega_i)} \\
&\leq\|\pi(\beta_{glo}^{(i)})\|^2_{s(\omega_i)} + \cfrac{1}{\Lambda} \|\beta_{glo}^{(i)}\|^2_{a(\omega_i)}.
\end{align*}
Therefore, we have
\[
\sum_{i=1}^{N_c} \Big(\|\beta_{glo}^{(i)}-\beta_{ms}^{(i)}\|_{a}^2 + \|\pi(\beta_{glo}^{(i)}-\beta_{ms}^{(i)})\|_{s}^2 \Big)\leq 2E(1+\Lambda^{-1})\sum_{i=1}^{N_c} \|u-u^{(m)}_{ms}\|^2_{a(\omega_i)}\leq 2ME(1+\Lambda^{-1})\|u-u^{(m)}_{ms}\|^2_{a}.
\]

\noindent
{\bf Step 4}:

In this final step, we will prove the required convergence.
Let $I = \{1,\cdots, k\}$.
From the adaptive enrichment algorithm, we add the online basis functions $\beta^{(i)}_{ms}$ for $i\in I$.
We will take $w$ in (\ref{eq:galerkin}) as
\begin{equation*}
w = u_{ms}^{(m)}-\sum_{i\in I}\beta_{ms}^{(i)}+\sum_{i=1}^{N} \sum_{j=1}^{J_i}c_{j}^{(i)}\psi_{j,ms}^{(i)} \in V_{ms}^{(m+1)}.
\end{equation*}
Using (\ref{eq:galerkin}) and (\ref{eq:rep}), we have
\begin{align*}
\|u-u_{ms}^{(m+1)}\|_{a}^{2} & \leq\|u-u_{ms}^{(m)}+\sum_{i\in I}\beta_{ms}^{(i)}-\sum_{i=1}^N\sum_{j=1}^{J_i}c_{j}^{(i)}\psi_{j,ms}^{(i)}\|_{a}^{2}\\
 & = \|\sum_{i\in I}(\beta_{glo}^{(i)}-\beta_{ms}^{(i)}) + \sum_{i\notin I}\beta_{glo}^{(i)} + \sum_{i=1}^N\sum_{j=1}^{J_i}c_{j}^{(i)}(\psi_{j}^{(i)} - \psi_{j,ms}^{(i)})\|_{a}^{2} \\
 & \leq 3\Big\{  \| \sum_{i\in I}(\beta_{glo}^{(i)}-\beta_{ms}^{(i)})\|_a^2
 + \| \sum_{i\notin I}\beta_{glo}^{(i)}\|_a^2
 + \| \sum_{i=1}^N\sum_{j=1}^{J_i}c_{j}^{(i)}(\psi_{j}^{(i)} - \psi_{j,ms}^{(i)}) \|_a^2 \Big\}.
\end{align*}
Using Step 2 and Step 3 as well as Lemma \ref{lem:on}, we see that
\begin{equation*}
 \| \sum_{i\in I}(\beta_{glo}^{(i)}-\beta_{ms}^{(i)})\|_a^2 +  \| \sum_{i=1}^N\sum_{j=1}^{J_i}c_{j}^{(i)}(\psi_{j}^{(i)} - \psi_{j,ms}^{(i)}) \|_a^2
 \leq C (\ell+1)^{d} E(1+\Lambda^{-1})\|u-u^{(m)}_{ms}\|^2_{a}.
\end{equation*}

Next, we will estimate the remaining term $\|\sum_{i\notin I} \beta_{glo}^{(i)}\|^2_a$. We write $p := \sum_{i\notin I} \beta_{glo}^{(i)}$. Then, by (\ref{eq:global_on_basis}),
\begin{align*}
\|p\|_{a}^{2} + \|\pi(p)\|_{s}^{2} &=   r(\sum_{i\notin I}\chi_i p)
 \leq \sum_{i\notin I}\Big(\sup_{v\in V_0(\omega_i)} \cfrac{r(v)}{\|v\|_a}\Big) \|\chi_i p\|_a \\
& \leq \sqrt{2}\sum_{i\notin I} \|z_i\|_{a^*}\Big(\|p\|^2_{a(\omega_i)}+\|p\|^2_{s(\omega_i)}\Big)^{\frac{1}{2}}\\
& \leq \sqrt{2} (1+\Lambda^{-1})^{\frac{1}{2}} \sum_{i\notin I}\|z_i\|_{a^*}\Big(\|p\|^2_{a(\omega_i)}+\|\pi(p)\|^2_{s(\omega_i)}\Big)^{\frac{1}{2}}.
\end{align*}
Thus, we have
\[
\|\sum_{i\notin I} \beta_{glo}^{(i)}\|^2_a \leq 2M(1+\Lambda^{-1})\sum_{i\notin I} \|z_i\|^2_{a^*}\leq 2M (1+\Lambda^{-1}) \theta \sum_{i=1}^{N_c} \|z_i\|^2_{a^*}
\]
where $\theta$ is defined in (\ref{eq:theta}).
Lastly, we will estimate $\|z_i\|_{a^*}$. By definition,
\[
z_i(v) = a(u-u^{(m)}_{ms},v)\leq \|u-u^{(m)}_{ms}\|_{a(\omega_i)}\|v\|_a \\
\]
Thus, we have
\begin{align*}
\sum_{i=1}^{N_c} \|z_i\|^2_{a^*} \leq   \sum_{i=1}^{N_c} \|u-u^{(m)}_{ms}\|^2_{a(\omega_i)}  \leq  M \|u-u^{(m)}_{ms}\|^2_a.
\end{align*}
Combining the above equations, we have
\begin{align*}
\|u-u_{ms}^{(m+1)}\|_{a}^{2} & \leq 3\Big(C(\ell+1)^{d}E(1+\Lambda^{-1})\|u-u_{ms}^{(m)}\|_{a}^{2} + 2M^2(1+\Lambda^{-1}) \theta \|u-u^{(m)}_{ms}\|^2_a\Big)\\
 & \leq 3(1+\Lambda^{-1})\Big(C(\ell+1)^{d}E + 2M^2\theta \Big)\|u-u^{(m)}_{ms}\|^2_a.
\end{align*}
This completes the proof.


\section{Conclusions}
\label{sec:conclusion}

In this paper, we develop an online adaptive enrichment algorithm for CEM-GMsFEM.
The CEM-GMsFEM, developed in \cite{chung2017constraint}, provides a general methodology for constructing
multiscale basis functions that give a mesh-dependent convergence rate, regardless the contrast and heterogeneties of the media.
In some applications, there is a need to further reduce the error without refining the coarse mesh. In these cases, one needs
to add basis functions.
We propose a strategy to compute new basis functions in the online stage using local residuals. These online basis functions
provide very fast (exponential) decay in error. Moreover, an adaptive strategy is proposed to enrich the basis
in some selected regions with large residuals. This strategy is determined by a user-defined parameter.
We show that this used-defined parameter relates directly to the convergence rate of the method.
That is, one can determine the convergence rate by using this parameter correspondingly.
Several numerical tests are shown to validate our estimates.

\bibliographystyle{plain}
\bibliography{references,references1,references_outline}

\begin{thebibliography}{10}

\bibitem{chan2016adaptive}
Ho~Yuen Chan, Eric Chung, and Yalchin Efendiev.
\newblock Adaptive mixed {GM}s{FEM} for flows in heterogeneous media.
\newblock {\em Numerical Mathematics: Theory, Methods and Applications},
  9(4):497--527, 2016.

\bibitem{chung2016adaptive}
Eric Chung, Yalchin Efendiev, and Thomas~Y Hou.
\newblock Adaptive multiscale model reduction with generalized multiscale
  finite element methods.
\newblock {\em Journal of Computational Physics}, 320:69--95, 2016.

\bibitem{chung2015residual}
Eric~T Chung, Yalchin Efendiev, and Wing~Tat Leung.
\newblock Residual-driven online generalized multiscale finite element methods.
\newblock {\em Journal of Computational Physics}, 302:176--190, 2015.

\bibitem{chung2017constraint}
Eric~T Chung, Yalchin Efendiev, and Wing~Tat Leung.
\newblock Constraint energy minimizing generalized multiscale finite element
  method.
\newblock {\em arXiv preprint arXiv:1704.03193}, 2017.

\bibitem{chung2016onlineperforated}
Eric~T Chung, Yalchin Efendiev, Wing~Tat Leung, Maria Vasilyeva, and Yating
  Wang.
\newblock Online adaptive local multiscale model reduction for heterogeneous
  problems in perforated domains.
\newblock {\em Applicable Analysis}, pages 1--30, 2016.

\bibitem{chung2015online}
ET~Chung, Y~Efendiev, and WT~Leung.
\newblock An online generalized multiscale discontinuous {Galerkin} method
  {(GMsDGM)} for flows in heterogeneous media.
\newblock {\em arXiv preprint arXiv:1504.04417}, 2015.

\bibitem{dur91}
L.J. Durlofsky.
\newblock Numerical calculation of equivalent grid block permeability tensors
  for heterogeneous porous media.
\newblock {\em Water Resour. Res.}, 27:699--708, 1991.

\bibitem{ee03}
W.~E and B.~Engquist.
\newblock Heterogeneous multiscale methods.
\newblock {\em Comm. Math. Sci.}, 1(1):87--132, 2003.

\bibitem{ge09_2}
Y.~Efendiev and J.~Galvis.
\newblock A domain decomposition preconditioner for multiscale high-contrast
  problems.
\newblock In Y.~Huang, R.~Kornhuber, O.~Widlund, and J.~Xu, editors, {\em
  Domain Decomposition Methods in Science and Engineering XIX}, volume~78 of
  {\em Lect. Notes in Comput. Science and Eng.}, pages 189--196.
  Springer-Verlag, 2011.

\bibitem{egh12}
Y.~Efendiev, J.~Galvis, and T.~Hou.
\newblock Generalized multiscale finite element methods.
\newblock {\em Journal of Computational Physics}, 251:116--135, 2013.

\bibitem{eh09}
Y.~Efendiev and T.~Hou.
\newblock {\em Multiscale Finite Element Methods: Theory and Applications}.
\newblock Springer, 2009.

\bibitem{ehw99}
Y.~Efendiev, T.~Hou, and X.H. Wu.
\newblock Convergence of a nonconforming multiscale finite element method.
\newblock {\em SIAM J. Numer. Anal.}, 37:888--910, 2000.

\bibitem{hw97}
T.~Hou and X.H. Wu.
\newblock A multiscale finite element method for elliptic problems in composite
  materials and porous media.
\newblock {\em J. Comput. Phys.}, 134:169--189, 1997.

\bibitem{hou2017}
Tom Hou and Pengchuan Zhang.
\newblock private communications.

\bibitem{hughes1995multiscale}
TJR Hughes.
\newblock Multiscale phenomena: {Green's functions, the Dirichlet-to-Neumann}
  formulation, subgrid scale models, bubbles and the origins of stabilized
  methods.
\newblock {\em Computer methods in applied mechanics and engineering},
  127(1):387--401, 1995.

\bibitem{maalqvist2014localization}
Axel M{\aa}lqvist and Daniel Peterseim.
\newblock Localization of elliptic multiscale problems.
\newblock {\em Mathematics of Computation}, 83(290):2583--2603, 2014.

\bibitem{owhadi2017multigrid}
Houman Owhadi.
\newblock Multigrid with rough coefficients and multiresolution operator
  decomposition from hierarchical information games.
\newblock {\em SIAM Review}, 59(1):99--149, 2017.

\bibitem{owhadi2014polyharmonic}
Houman Owhadi, Lei Zhang, and Leonid Berlyand.
\newblock Polyharmonic homogenization, rough polyharmonic splines and sparse
  super-localization.
\newblock {\em ESAIM: Mathematical Modelling and Numerical Analysis},
  48(2):517--552, 2014.

\bibitem{Wheeler_mortar_MS_12}
M.F. Wheeler, G.~Xue, and I.~Yotov.
\newblock A multiscale mortar multipoint flux mixed finite element method.
\newblock {\em ESAIM Math. Model. Numer. Anal.}, 46(4):759--796, 2012.

\bibitem{weh02}
X.H. Wu, Y.~Efendiev, and T.Y. Hou.
\newblock Analysis of upscaling absolute permeability.
\newblock {\em Discrete and Continuous Dynamical Systems, Series B.},
  2:158--204, 2002.

\end{thebibliography}

\appendix

\section{Proof of Lemma \ref{lem:sum}}

In this appendix, we give a proof for Lemma \ref{lem:sum}.
Let $K_i$ be a coarse element and let $K_{i,n}$ be the oversampled region by enlarging $K_{i}$ by $n$ coarse grid layers.
For integers $m>n$, we denote $ \chi^{m,n}_i$ as the cutoff function used in \cite{chung2017constraint}.
In particular, this function satisfies $\chi^{m,n}_i=1$ in $K_{i,n}$ and $\chi^{m,n}_i=0$ in $\Omega \backslash K_{i,m}$.
Next, we set $w^{(i)}$ as $\sum_{j=1}^{J_i}c^{(i)}_j(\psi_j^{(i)} - \psi_{j,ms}^{(i)})$ and $w$ as $\sum_{i=1}^N w^{(i)}$. Then, by the definitions (\ref{eq:local_basis}) and (\ref{eq:global_basis}),
we have
\[
a((1-\chi^{\ell+1,\ell}_i)w,\sum_{j=1}^{J_i}c^{(i)}_j \psi^{(i)}_{j,ms}) + s(\pi((1-\chi^{\ell+1,\ell}_i)w),\pi(\sum_{j=1}^{J_i}c^{(i)}_j \psi^{(i)}_{j,ms}))=0
 \]
 and
 \[
 a((1-\chi^{\ell+1,\ell}_i)w,\sum_{j=1}^{J_i}c^{(i)}_j \psi^{(i)}_{j}) + s(\pi((1-\chi^{\ell+1,\ell}_i)w),\pi(\sum_{j=1}^{J_i}c^{(i)}_j \psi^{(i)}_{j}))= s(\sum_{j=1}^{J_i}c^{(i)}_j \phi^{(i)}_{j},(1-\chi^{\ell+1,\ell}_i)w) = 0
 \]
  since $\text{supp}\{\sum_{j=1}^{J_i}c^{(i)}_j \phi^{(i)}_{j,ms}\}\subset \text{supp}\{\sum_{j=1}^{J_i}c^{(i)}_j \psi^{(i)}_{j,ms}\}\subset K_{i,\ell}$.
Therefore, by subtracting the above two equations,
we have
\begin{align*}
\|w\|^2_a +\|\pi(w)\|^2_s &= \sum_{i=1}^N a((\chi^{\ell+1,\ell}_i)w,w^{(i)}) + s(\pi((\chi^{\ell+1,\ell}_i)w),\pi(w^{(i)}))\\
& \leq \sum_{i=1}^N \Big(\|(\chi^{\ell+1,\ell}_i)w\|_a \|w^{(i)}\|_a + \|\pi((\chi^{\ell+1,\ell}_i)w)\|_s \|\pi(w^{(i)})\|_s\Big).
\end{align*}
Now, we estimate the right hand side of the above.
Notice that,
\begin{align*}
\|(\chi^{\ell+1,\ell}_i)w\|^2_a + \|\pi((\chi^{\ell+1,\ell}_i)w)\|^2_s &\leq C(\|w\|^2_{a(K_{i,\ell+1})}+\|w\|^2_{s(K_{i,\ell+1})})\\
&\leq C(1+\Lambda^{-1})(\|w\|^2_{a(K_{i,\ell+1})}+\|\pi(w)\|^2_{s(K_{i,\ell+1})}).
\end{align*}
Summing over $i=1,2,\cdots, N$,
we have
 \begin{align*}
\sum_{i=1}^N \|(\chi^{\ell+1,\ell}_i)w\|^2_a + \|\pi((\chi^{\ell+1,\ell}_i)w)\|^2_s &\leq C(1+\Lambda^{-1})\sum_{i=1}^N (\|w\|^2_{a(K_{i,\ell+1})}+\|\pi(w)\|^2_{s(K_{i,\ell+1})}) \\
& \leq C(1+\Lambda^{-1})(\ell+1)^d(\|w\|^2_{a}+\|\pi(w)\|^2_{s}).
\end{align*}
Hence, we conclude that
\[
\|w\|^2_a +\|\pi(w)\|^2_s \leq C(1+\Lambda^{-1})(\ell+1)^d \sum_{i=1}^N \Big(\|w^{(i)}\|_a^2 + \|\pi(w^{(i)})\|_s^2\Big).
\]
This completes the proof of the lemma.

\end{document}